\documentclass[11pt,a4paper]{amsart}
\usepackage[T1]{fontenc}
\usepackage[utf8]{inputenc}
\usepackage{lmodern}
\usepackage{amsmath,amssymb,mathtools}
\usepackage{microtype}
\usepackage{float}
\usepackage{pgfplots}
\pgfplotsset{compat=1.18}
\usepackage[hidelinks]{hyperref}
\numberwithin{equation}{section}
\theoremstyle{plain}
\newtheorem{theorem}{Theorem}[section]
\newtheorem{lemma}[theorem]{Lemma}
\theoremstyle{remark}
\newtheorem{remark}[theorem]{Remark}
\newtheorem{example}[theorem]{Example}

\newcommand{\cft}[2]{[z^{#1}]#2}
\expandafter\def\csname plotk1\endcsname{(1,3.50000000) (2,3.75000000) (3,3.87500000) (4,3.93750000) (5,3.96875000) (6,3.98437500) (7,3.99218750) (8,3.99609375) (9,3.99804688) (10,3.99902344) (11,3.99951172) (12,3.99975586)}
\expandafter\def\csname plotk2\endcsname{(1,5.50000000) (2,6.75000000) (3,7.37500000) (4,7.68750000) (5,7.84375000) (6,7.92187500) (7,7.96093750) (8,7.98046875) (9,7.99023438) (10,7.99511719) (11,7.99755859) (12,7.99877930)}
\expandafter\def\csname plotk3\endcsname{(1,8.00000000) (2,11.75000000) (3,13.87500000) (4,14.93750000) (5,15.46875000) (6,15.73437500) (7,15.86718750) (8,15.93359375) (9,15.96679688) (10,15.98339844) (11,15.99169922) (12,15.99584961)}
\expandafter\def\csname plotk4\endcsname{(1,11.00000000) (2,20.00000000) (3,25.87500000) (4,28.93750000) (5,30.46875000) (6,31.23437500) (7,31.61718750) (8,31.80859375) (9,31.90429688) (10,31.95214844) (11,31.97607422) (12,31.98803711)}
\expandafter\def\csname plotk6\endcsname{(1,18.50000000) (2,55.50000000) (3,88.75000000) (4,108.00000000) (5,117.96875000) (6,122.98437500) (7,125.49218750) (8,126.74609375) (9,127.37304688) (10,127.68652344) (11,127.84326172) (12,127.92163086)}
\expandafter\def\csname plotN20\endcsname{(1,0) (2,0) (3,0) (4,1) (5,3) (6,7) (7,14) (8,7) (9,7) (10,0) (11,13) (12,59)}
\expandafter\def\csname plotN21\endcsname{(1,0) (2,0) (3,0) (4,1) (5,3) (6,7) (7,14) (8,6) (9,5) (10,-4) (11,7) (12,51)}

\title[Brown-Gap Healing and Completeness Thresholds]{Brown-Gap Healing and Exact Completeness Thresholds\\for Leading-Unit Positive Linear Recurrences}
\author{Tony Newton}
\address{Newton Astro Labs (independent research), United Kingdom}
\email{tony.newton79@gmail.com}
\subjclass[2020]{11B37, 11B39, 11B83}
\keywords{positive linear recurrence sequences, complete sequences, Brown's criterion, generating functions, boundary equality, monotonicity}

\begin{document}
\begin{abstract}
Consider the standard positive linear recurrence sequence with $g$ leading unit coefficients, $k$ subsequent zeros and final coefficient $N$, where $g,k,N\ge 1$. We determine its exact completeness threshold. Put
\[
Q_g(z)=1-z-\cdots-z^g,\qquad
 a_m^{(g)}=\cft{m}{\frac{1}{(1-z)^2Q_g(z)}}.
\]
Then the sequence is complete precisely when
$1\le N\le \lfloor (a_{g+k}^{(g)}+1)/2^g\rfloor$.
The proof analyzes Brown's gaps and establishes their nonnegativity beyond a finite initial block, including when the critical gap is zero. The unrounded threshold is strictly increasing in $g$, establishing leading-unit monotonicity of completeness.
\end{abstract}
\maketitle
\thispagestyle{plain}

\section{Introduction and main results}
A sequence of positive integers is \emph{complete} if every positive integer is a sum of distinct indexed terms of the sequence. Representations need not be unique. Brown's criterion~\cite{Brown1961} states that a nondecreasing positive sequence with $H_1=1$ is complete if and only if
\begin{equation}\label{eq:brown}
 B_n:=1+\sum_{i=1}^{n-1}H_i-H_n\ge0\qquad(n\ge1).
\end{equation}
Indeed, if the preceding inequalities hold, the first $n-1$ terms represent every integer from $0$ through $\sum_{i<n}H_i$; the term $H_n$ extends this interval without a gap precisely when $B_n\ge0$.

Fix integers $g,k,N\ge1$ and let $L=g+k+1$. We consider the standard positive linear recurrence sequence associated with the coefficient list
$[\underbrace{1,\ldots,1}_{g},\underbrace{0,\ldots,0}_{k},N]$:
\begin{equation}\label{eq:initial}
 H_1=1,\qquad H_n=1+\sum_{j=1}^{\min(g,n-1)}H_{n-j}\quad(2\le n\le L),
\end{equation}
\begin{equation}\label{eq:recurrence}
 H_n=\sum_{j=1}^{g}H_{n-j}+NH_{n-L}\qquad(n\ge L+1).
\end{equation}
These terms are strictly increasing, so Brown's criterion applies. Introduce the formal power series
\begin{equation}\label{eq:A}
 Q_g(z)=1-\sum_{j=1}^{g}z^j,\qquad
 A_g(z)=\frac1{(1-z)^2Q_g(z)}=\sum_{m\ge0}a_m^{(g)}z^m,
\end{equation}
and the real threshold
\begin{equation}\label{eq:M}
 M_{g,k}=\frac{a_{g+k}^{(g)}+1}{2^g}.
\end{equation}
\begin{theorem}[Exact threshold]\label{thm:threshold}
For every $g,k,N\ge1$, the sequence \eqref{eq:initial}--\eqref{eq:recurrence} is complete if and only if
\begin{equation}\label{eq:threshold}
 1\le N\le \lfloor M_{g,k}\rfloor
 =\left\lfloor\frac{1+\cft{g+k}{((1-z)^2Q_g(z))^{-1}}}{2^g}\right\rfloor.
\end{equation}
If $M_{g,k}$ is an integer, equality $N=M_{g,k}$ is included.
\end{theorem}
\begin{theorem}[Leading-unit monotonicity]\label{thm:monotonicity}
For all $g,k\ge1$, one has $M_{g+1,k}>M_{g,k}$. Consequently, if $[1^g,0^k,N]$ is complete, then $[1^{g+1},0^k,N]$ is complete.
\end{theorem}

The one-leading-unit classification occurs in~\cite[Theorem~1.15]{Boldyriew2020manuscript}. A proof of the two-leading-unit Fibonacci threshold appears in the September 2020 manuscript~\cite[Theorem~1.16]{Boldyriew2020manuscript}, while the later arXiv revision~\cite[Conjecture~3.1]{Boldyriew2021arxiv} labels that same threshold a conjecture. The implication in Theorem~\ref{thm:monotonicity} was posed in~\cite[Conjecture~3.2]{Boldyriew2021arxiv}; see also the published account~\cite{Boldyriew2023} and introductory survey~\cite{Boldyriew2020survey}. The contribution here is a uniform coefficient criterion for arbitrary $g$ and $k$, including the equality case, together with a direct proof of strict monotonicity for the real threshold.

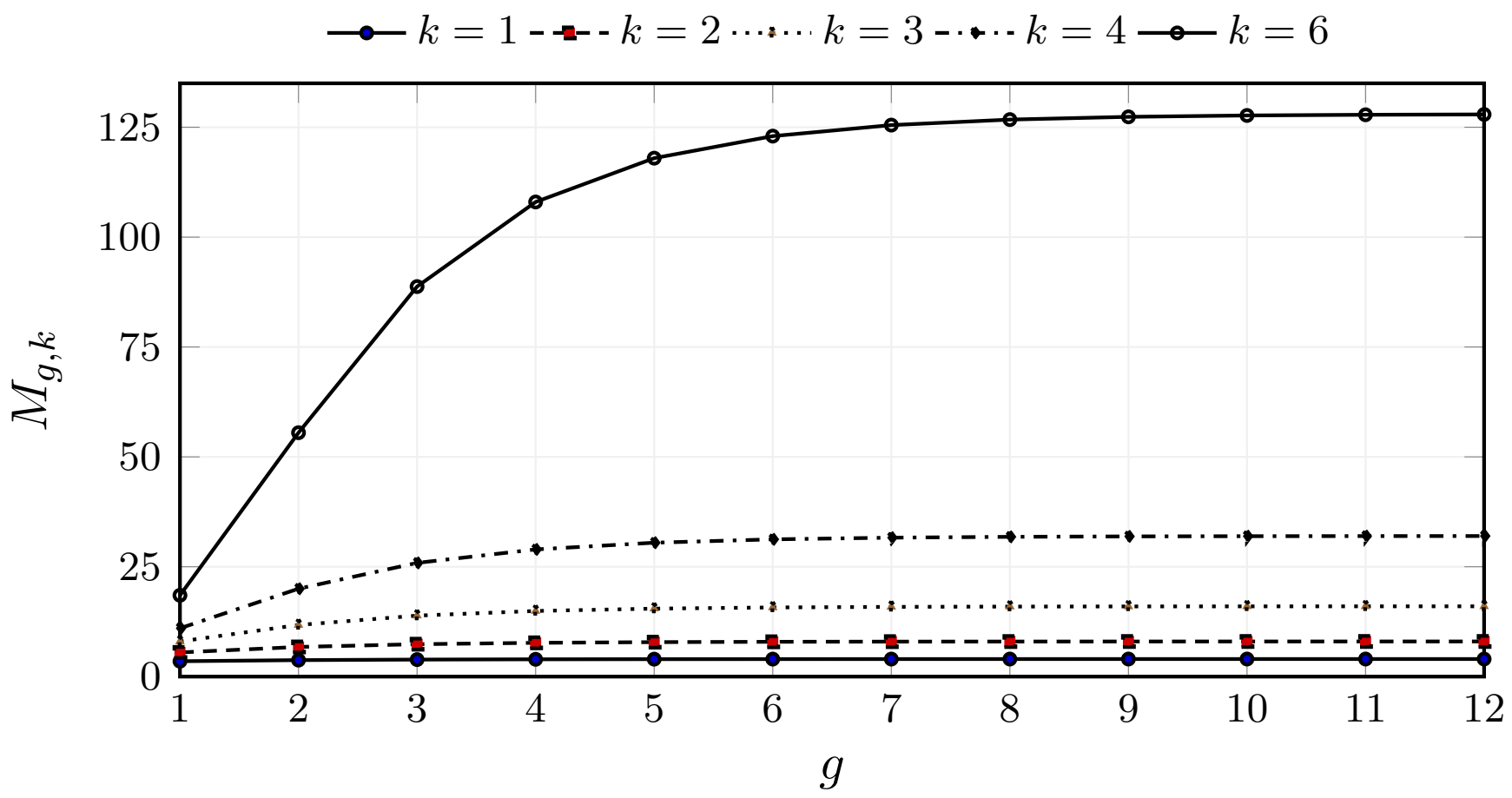
\begin{figure}[H]
\centering
\begin{tikzpicture}
\begin{axis}[
 width=.91\linewidth,height=6.1cm,
 xlabel={$g$},ylabel={$M_{g,k}$},xmin=1,xmax=12,ymin=0,ymax=135,
 xtick={1,...,12},ytick={0,25,50,75,100,125},
 legend style={at={(0.5,1.02)},anchor=south,legend columns=5,draw=none,font=\footnotesize},
 tick label style={font=\footnotesize},line width=.8pt,
 grid=major,grid style={gray!12},mark size=1.2pt]
\addplot+[mark=*,black] coordinates {\csname plotk1\endcsname};\addlegendentry{$k=1$}
\addplot+[mark=square*,black,densely dashed] coordinates {\csname plotk2\endcsname};\addlegendentry{$k=2$}
\addplot+[mark=triangle*,black,dotted] coordinates {\csname plotk3\endcsname};\addlegendentry{$k=3$}
\addplot+[mark=diamond*,black,dashdotted] coordinates {\csname plotk4\endcsname};\addlegendentry{$k=4$}
\addplot+[mark=o,black,thick] coordinates {\csname plotk6\endcsname};\addlegendentry{$k=6$}
\end{axis}
\end{tikzpicture}
\caption{The unrounded completeness thresholds $M_{g,k}$ for fixed $k$ and $1\le g\le12$. The strict growth underlying Theorem~\ref{thm:monotonicity} is visible before taking integer parts.}
\label{fig:thresholds}
\end{figure}

\section{The Brown-gap generating function}
For fixed $g$, write $a_m=a_m^{(g)}$ and set $a_m=0$ for $m<0$. The identity $(1-z)Q_g(z)=1-2z+z^{g+1}$ gives
\begin{equation}\label{eq:a-rec}
 a_m-2a_{m-1}+a_{m-g-1}=1\qquad(m\ge0),
\end{equation}
and hence
\begin{equation}\label{eq:a-init}
 a_r=2^{r+1}-1\qquad(0\le r\le g).
\end{equation}
Every $a_m$ is a positive integer because $(1-z)^{-2}$ and $Q_g(z)^{-1}$ both have nonnegative integer coefficients and positive constant term. For $0\le r\le g$, put $R_r=(a_{k+r}+1)/2^r$. Then \eqref{eq:a-rec} yields
\begin{equation}\label{eq:R-diff}
 R_{r-1}-R_r=\frac{a_{k+r-g-1}}{2^r}\ge0\qquad(1\le r\le g).
\end{equation}
The final difference is strict because $a_{k-1}>0$. Therefore
\begin{equation}\label{eq:R-chain}
 R_0\ge R_1\ge\cdots\ge R_g=M_{g,k}<a_k+1.
\end{equation}
Write $D(z)=Q_g(z)-Nz^L$. The initial conditions and recurrence give
\[
 H(z):=\sum_{n\ge1}H_nz^{n-1}=\frac{1-z^L}{(1-z)D(z)}.
\]
Consequently, the generating function of the gaps in \eqref{eq:brown} is
\begin{align}\label{eq:B}
 B(z)&:=\sum_{n\ge1} B_nz^{n-1}
 =\frac1{1-z}+\frac{2z-1}{1-z}H(z)\\
 &\notag=z^{g+1}F(z),
\end{align}
where
\begin{equation}\label{eq:F}
 F(z)=\frac{1-(N-1)z^k+(N-2)z^{k+1}}{(1-z)^2D(z)}
 =\sum_{m\ge0} f_mz^m.
\end{equation}
Set $f_m=0$ for $m<0$. Equation \eqref{eq:B} gives
\begin{equation}\label{eq:shift}
 B_n=0\quad(1\le n\le g+1),\qquad B_{m+g+2}=f_m\quad(m\ge0).
\end{equation}
\begin{lemma}[Initial block]\label{lem:initial}
For $0\le m<k$, one has $f_m=a_m$. For $0\le r\le g$,
\begin{equation}\label{eq:initial-block}
 f_{k+r}=a_{k+r}+1-N2^r.
\end{equation}
Thus $N\le M_{g,k}$ implies $f_m\ge0$ for $m<L$. Conversely, completeness implies $N\le M_{g,k}$.
\end{lemma}
\begin{proof}
Up to degree $L-1$, replacing $D(z)^{-1}$ in \eqref{eq:F} by $Q_g(z)^{-1}$ does not change any coefficient. This gives $f_m=a_m$ for $m<k$ and, for $0\le r\le g$,
\[
 f_{k+r}=a_{k+r}-(N-1)a_r+(N-2)a_{r-1}.
\]
Substitution of \eqref{eq:a-init} proves \eqref{eq:initial-block}. By \eqref{eq:R-chain}, $N\le M_{g,k}\le R_r$ makes the entire block nonnegative. Conversely, the gap at the critical index is
$B_{2g+k+2}=f_{L-1}=a_{g+k}+1-N2^g$; Brown's criterion requires it to be nonnegative.
\end{proof}
Multiplication of \eqref{eq:F} by $D(z)$ and coefficient extraction give the inhomogeneous recurrence
\begin{equation}\label{eq:f-rec}
 f_m=\sum_{j=1}^{g}f_{m-j}+Nf_{m-L}+c,\qquad c=k+2-N,\qquad m\ge k.
\end{equation}
Its constant forcing term makes it possible to control all gaps after the initial block.

\section{Healing of the gaps with at least two leading units}
\begin{lemma}\label{lem:g-ge2}
Suppose $g\ge2$, $k\ge1$ and $1\le N\le M_{g,k}$. Then $f_m\ge0$ for all $m$, and $f_m>0$ whenever $m\ge L$.
\end{lemma}
\begin{proof}
Lemma~\ref{lem:initial} supplies $f_m\ge0$ for $m<L$. If $k=1$, \eqref{eq:R-chain} and $a_1=3$ imply $N\le M_{g,1}<4$, so $N\le3$ and $c=3-N\ge0$. Moreover, $M_{g,1}=4-2^{-g}$ is not an integer; thus \eqref{eq:initial-block} shows that the whole initial block is strictly positive. Recurrence \eqref{eq:f-rec} preserves nonnegativity and, because the delayed term $Nf_{m-L}$ is positive at $m=L$ and thereafter, strict positivity follows by induction.

Suppose $k\ge2$. Since $g\ge2$,
\begin{equation}\label{eq:seed}
 f_0=1,\qquad f_1=3,\qquad f_2\ge1.
\end{equation}
Indeed, $f_2=a_2=7$ for $k\ge3$, while for $k=2$, $f_2=a_2+1-N=8-N\ge1$ by \eqref{eq:R-chain}. Using \eqref{eq:f-rec} and the nonnegative initial block,
\begin{align}
 f_L&\ge Nf_0+c=k+2,\label{eq:f-L}\\
 f_{L+1}&\ge f_L+Nf_1+c=f_L+2N+k+2,\label{eq:f-L1}\\
 f_{L+2}&\ge f_{L+1}+f_L+Nf_2+c\ge f_{L+1}+f_L+k+2.\label{eq:f-L2}
\end{align}
Let $T=\max\{0,N-k-2\}$, so that $c\ge-T$. In particular, the last two quantities in \eqref{eq:f-L1}--\eqref{eq:f-L2} exceed $T$. If $m\ge L+3$ and all earlier coefficients are nonnegative with $f_{m-1},f_{m-2}\ge T$, then
\[
 f_m\ge f_{m-1}+f_{m-2}+c\ge T.
\]
Induction proves nonnegativity. Strict positivity is immediate when $T>0$. When $T=0$, one has $c\ge0$ and the two positive preceding coefficients force $f_m>0$ throughout the tail.
\end{proof}

\section{The one-leading-unit boundary}
The preceding argument uses two consecutive positive terms, so the case $g=1$ requires a different induction.
\begin{lemma}\label{lem:g1}
If $g=1$, $k\ge1$ and $1\le N\le M_{1,k}$, then $f_m\ge0$ for all $m$, and $f_m>0$ for $m\ge L$. This includes integer equality $N=M_{1,k}$.
\end{lemma}
\begin{proof}
Here $L=k+2$, $Q_1(z)=1-z$, and
\begin{equation}\label{eq:g1}
 a_m=\binom{m+2}{2},\qquad M_{1,k}=\frac{L(L+1)+2}{4}.
\end{equation}
By Lemma~\ref{lem:initial}, $f_m\ge0$ for $m<L$. In fact $f_m\ge1$ for $m\le L-2$: this follows from $f_m=a_m>0$ for $m<k$, while for $m=k=L-2$ it follows from $N\le M_{1,k}<a_k+1$ and integrality. Only $f_{L-1}$ may vanish.

The recurrence is now
\begin{equation}\label{eq:g1rec}
 f_m=f_{m-1}+Nf_{m-L}+L-N\qquad(m\ge k).
\end{equation}
For $L\le m\le2L-2$, the delayed coefficient is at least one, so $f_m\ge f_{m-1}+L$; hence $f_{2L-2}\ge(L-1)L$. At $m=2L-1$, the delayed coefficient is $f_{L-1}$ and may be zero, but
\[
 f_{2L-1}\ge f_{2L-2}+L-N\ge L^2-N>0,
\]
because $N\le M_{1,k}<L^2$ for $L\ge3$. For $m\ge2L$ the delayed coefficient has index at least $L$ and is positive by induction; in \eqref{eq:g1rec} this gives $f_m\ge f_{m-1}+L>0$.
\end{proof}
\begin{proof}[Proof of Theorem~\ref{thm:threshold}]
Necessity follows from Lemma~\ref{lem:initial}. If $N\le M_{g,k}$, that lemma supplies the initial nonnegative coefficients. Lemma~\ref{lem:g-ge2} for $g\ge2$, or Lemma~\ref{lem:g1} for $g=1$, supplies all remaining coefficients. Thus every Brown gap in \eqref{eq:shift} is nonnegative. Brown's criterion proves completeness, and integrality of $N$ gives \eqref{eq:threshold}.
\end{proof}

\section{Strict monotonicity of the threshold}
\begin{lemma}\label{lem:monotonicity}
$M_{g+1,k}>M_{g,k}$ for every $g,k\ge1$.
\end{lemma}
\begin{proof}
Write $A=A_g$, $\widetilde A=A_{g+1}$, and $a_m=\cft{m}{A}$, $\widetilde a_m=\cft{m}{\widetilde A}$. Since $Q_{g+1}=Q_g-z^{g+1}$,
\begin{equation}\label{eq:Adiff}
 \widetilde A-A=\frac{z^{g+1}}{(1-z)^2Q_gQ_{g+1}}=z^{g+1}\frac{\widetilde A}{Q_g}.
\end{equation}
Put $n=g+k$. The $g+1$ version of \eqref{eq:a-rec} reads $\widetilde a_{n+1}-2\widetilde a_n+\widetilde a_{n-g-1}=1$. It follows that
\begin{align}\label{eq:monotonicity-id}
 (\widetilde a_{n+1}+1)-2(a_n+1)
 &=2(\widetilde a_n-a_n)-\widetilde a_{n-g-1}\notag\\
 &=\cft{k-1}{\widetilde A(z)\!\left(\frac{2}{Q_g(z)}-1\right)}>0.
\end{align}
The last inequality holds because $\widetilde A$ has strictly positive coefficients, while $2/Q_g-1$ has nonnegative coefficients and constant coefficient one. Dividing \eqref{eq:monotonicity-id} by $2^{g+1}$ yields $M_{g+1,k}-M_{g,k}>0$.
\end{proof}
\begin{proof}[Proof of Theorem~\ref{thm:monotonicity}]
The strict inequality is Lemma~\ref{lem:monotonicity}. If $[1^g,0^k,N]$ is complete, Theorem~\ref{thm:threshold} gives $N\le M_{g,k}<M_{g+1,k}$. The same theorem then proves completeness of $[1^{g+1},0^k,N]$.
\end{proof}
\begin{remark}
Although $M_{g,k}$ increases strictly, the integer thresholds $\lfloor M_{g,k}\rfloor$ need not increase at every step.
\end{remark}

\section{Special cases and an equality example}
When $g=1$, \eqref{eq:g1} gives the known threshold
\[
 \lfloor M_{1,k}\rfloor=\left\lceil\frac{(k+2)(k+3)}{4}\right\rceil.
\]
For $g=2$, define $F_1=1$, $F_2=2$ and $F_{m+1}=F_m+F_{m-1}$. Formal-series division yields $a_m^{(2)}=F_{m+4}-m-4$ and hence
\begin{equation}\label{eq:fib}
 \lfloor M_{2,k}\rfloor=\left\lfloor\frac{F_{k+6}-k-5}{4}\right\rfloor.
\end{equation}
For $1\le k\le g$, iteration of \eqref{eq:a-rec} from \eqref{eq:a-init} gives
\[
 a_{g+k}^{(g)}=2^{g+k+1}-(k-1)2^k-2,\qquad
 \lfloor M_{g,k}\rfloor=2^{k+1}-\left\lceil\frac{k}{2^{g-k}}\right\rceil,
\]
which agrees with the range already covered in~\cite[Theorem~1.17]{Boldyriew2020manuscript}. These are consistency checks on the general formula.

\begin{example}[Vanishing critical gap]\label{ex:equality}
Take $(g,k,N)=(2,4,20)$, so $L=7$ and $M_{2,4}=20$. The sequence begins
\[
 1,\ 2,\ 4,\ 7,\ 12,\ 20,\ 33,\ 73,\ 146,\ 299,\ 585,\ 1124,\ldots .
\]
At the critical index, $B_{10}=0$, but the next gaps are $B_{11}=13$, $B_{12}=59$, and $B_{13}=198$. For $N=21$ with the same $g,k$, the critical gap is $B_{10}=-4$, so the corresponding sequence is incomplete. Thus equality in the threshold can be admissible even when the first decisive gap vanishes.
\end{example}

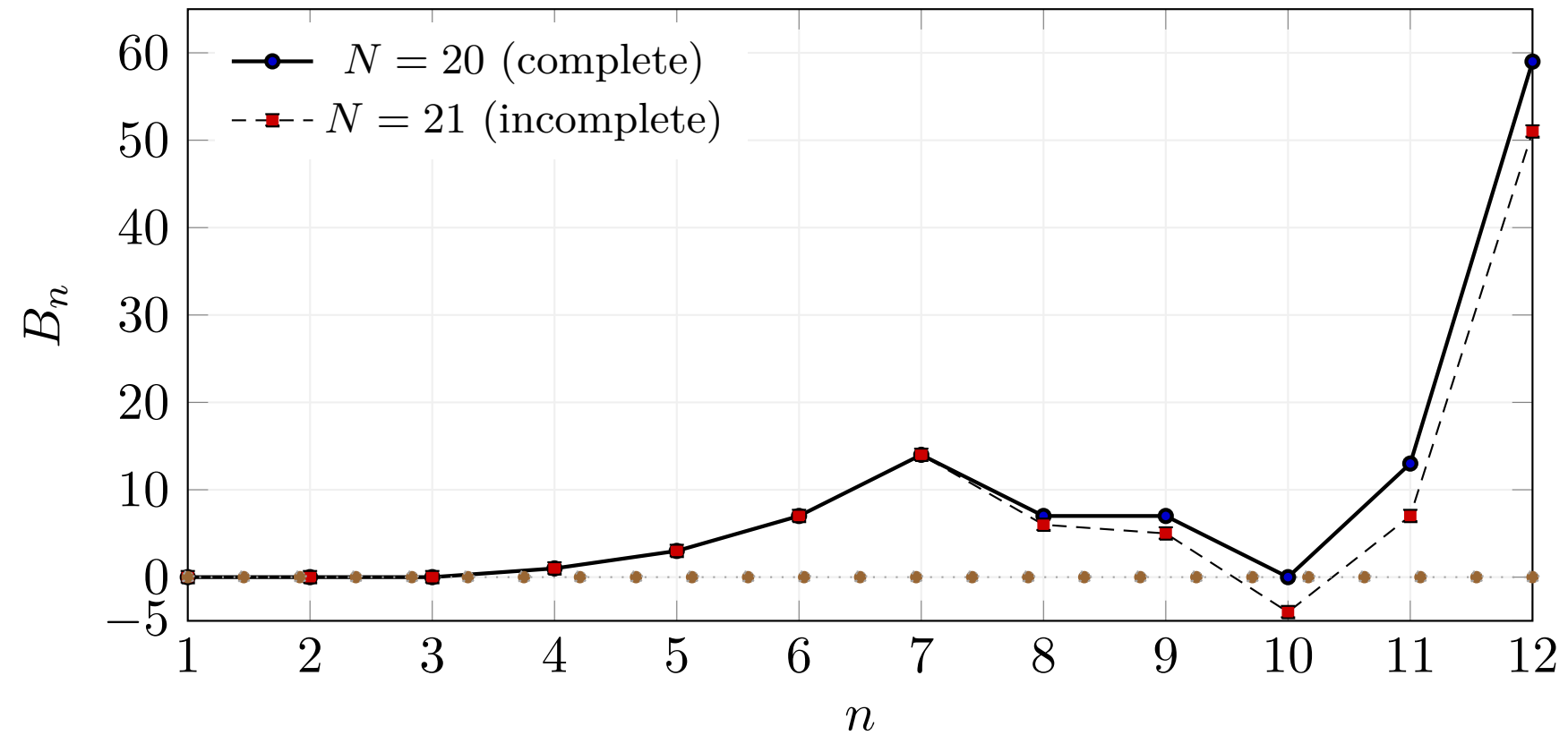
\begin{figure}[H]
\centering
\begin{tikzpicture}
\begin{axis}[
 width=.91\linewidth,height=6.1cm,
 xlabel={$n$},ylabel={$B_n$},xmin=1,xmax=12,ymin=-5,ymax=65,
 xtick={1,...,12},ytick={-5,0,10,20,30,40,50,60},
 legend style={at={(.02,.98)},anchor=north west,draw=none,font=\footnotesize},
 grid=major,grid style={gray!12},mark size=1.3pt,unbounded coords=jump]
\addplot+[black,thick,mark=*] coordinates {\csname plotN20\endcsname};\addlegendentry{$N=20$ (complete)}
\addplot+[black,densely dashed,mark=square*] coordinates {\csname plotN21\endcsname};\addlegendentry{$N=21$ (incomplete)}
\addplot+[black!35,dotted,forget plot,domain=1:12] {0};
\end{axis}
\end{tikzpicture}
\caption{Brown gaps for $(g,k)=(2,4)$ and two neighboring last coefficients. At $N=20$, $B_{10}=0$ and subsequent gaps are positive; at $N=21$, the critical gap is negative. Values are computed directly from \eqref{eq:initial}--\eqref{eq:recurrence}.}
\label{fig:boundary}
\end{figure}

\appendix
\section{Precise Mathematical Formulation of the Main Result}
Fix $g,k,N\in\mathbb{Z}_{\ge1}$ and $L=g+k+1$. Let $(H_n)_{n\ge1}$ be the uniquely determined integer sequence given by \eqref{eq:initial} and \eqref{eq:recurrence}. Completeness means that for every $t\in\mathbb{Z}_{\ge1}$ there is a finite set $S\subset\mathbb{Z}_{\ge1}$ such that $t=\sum_{i\in S}H_i$; repeated numerical values, if present, are distinguished by index. Put
\[
 Q_g(z)=1-\sum_{j=1}^{g}z^j,\qquad
 C(g,k)=\left\lfloor 2^{-g}\!\left(1+\cft{g+k}{\frac1{(1-z)^2Q_g(z)}}\right)\right\rfloor.
\]
All coefficient extractions are in the ring of formal power series $\mathbb Z[[z]]$. The principal assertion is the equivalence
\[
 (H_n)_{n\ge1}\ \text{is complete}
 \quad\Longleftrightarrow\quad 1\le N\le C(g,k).
\]
No strict inequality is imposed at the upper endpoint. For each fixed $k\ge1$, the real quantity $M_{g,k}$ in \eqref{eq:M} is strictly increasing as a function of the positive integer $g$. The proofs are given in Sections~2--5.

\section*{Acknowledgment and computational tools}
Computational algorithms, proof-search procedures, and certificate criteria were developed by the author/s. Generative AI systems, including ChatGPT and Claude, served as computational assistants and were used for debugging and cross-checking. Mathematical claims were accepted strictly on the basis of independently reproducible exact, symbolic, rational, interval, exhaustive, or controlled high-precision certificates, and were independently verified.

\end{document}